\documentclass[10pt]{article}
\usepackage{cite}
\usepackage{mathrsfs}
\usepackage{amsfonts}
\usepackage{amsmath}
\usepackage{amsfonts,amssymb}
\usepackage{dsfont}
\usepackage{curves}
\usepackage{mathrsfs}
\usepackage{pifont}
\usepackage{amssymb}
\allowdisplaybreaks

\numberwithin{equation}{section}

\date{}

\def\BigRoman{\uppercase\expandafter{\romannumeral\number\count 255 }}
\def\Romannumeral{\afterassignment\BigRoman\count255=}

\begin{document}
\title{Spanning trees and signless Laplacian spectral radius in graphs
%\thanks{Supported by }
}
\author{\small Sufang Wang$^{1}$, Wei Zhang$^{2}$\footnote{Corresponding
author. E-mail address: wangsufangjust@163.com (S. Wang), zw\_wzu@163.com (W. Zhang).}\\
\small $1$.  School of Public Management, Jiangsu University of Science and Technology,\\
\small Zhenjiang, Jiangsu 212100, China\\
\small $2$.  School of Economics and Management, Wenzhou University of Technology,\\
\small Wenzhou, Zhejiang 325000, China
}

\maketitle
\begin{abstract}
\noindent Let $G$ be a connected graph and let $k$ be a positive integer. Let $T$ be a spanning tree of $G$. The leaf degree of a vertex $v\in V(T)$
is defined as the number of leaves adjacent to $v$ in $T$. The leaf degree of $T$ is the maximum leaf degree among all the vertices of $T$. Let $A(G)$
be the adjacency matrix of $G$ and $D(G)$ be the diagonal degree matrix of $G$. Let $Q(G)=D(G)+A(G)$ be the signless Laplacian matrix of $G$. The
largest eigenvalue of $Q(G)$, denoted by $q(G)$, is called the signless Laplacian spectral radius of $G$. In this paper, we investigate the connection
between the spanning tree and the signless Laplacian spectral radius of $G$, and put forward a sufficient condition based upon the signless Laplacian
spectral radius to guarantee that a graph $G$ contains a spanning tree with leaf degree at most $k$. Finally, we construct some extremal graphs to claim
all the bounds obtained in this paper are sharp.
\\
\begin{flushleft}

(2020) Mathematics Subject Classification: 05C50, 05C05, 90B99
\end{flushleft}
\end{abstract}

\section{Introduction}

Graphs discussed in this paper are simple and undirected. Let $G$ be a graph with vertex set $V(G)=\{v_1,v_2,\ldots,v_n\}$ and edge set $E(G)$.
The order of $G$ is denoted by $|V(G)|=n$. A graph $G$ is called trivial if $n=1$.  The diagonal degree matrix of $G$ is denoted by
$D(G)=diag(d_1,d_2,\ldots,d_n)$, where $d_i$ denotes the degree of a vertex $v_i$ in $G$. The adjacency matrix of $G$ is defined to be the
$(0,1)$-matrix $A(G)=(a_{ij})_{n\times n}$, where $a_{ij}=1$ if and only if $v_i$ and $v_j$ are adjacent in $G$. Let $Q(G)=D(G)+A(G)$ be the
signless Laplacian matrix of $G$. The largest eigenvalues of $A(G)$ and $Q(G)$, denoted by $\rho(G)$ and $q(G)$, are called the spectral radius
and the signless Laplacian spectral radius of $G$, respectively.

For a vertex subset $S$ of $V(G)$, let $G[S]$ denote the subgraph of $G$ induced by $S$, and let $G-S$ be the subgraph derived from $G$ by
deleting the vertices in $S$ together with their incident edges. Let $G_1$ and $G_2$ be two vertex disjoint graphs. We denote by $G_1\cup G_2$
the disjoint union of $G_1$ and $G_2$, and by $G_1\vee G_2$ the join of $G_1$ and $G_2$. Let $K_n$ denote the complete graph of order $n$.

For two integers $a$ and $b$ with $1\leq a\leq b$, a spanning subgraph $F$ of $G$ is called an $[a,b]$-factor of $G$ if $a\leq d_F(v)\leq b$
for every $v\in V(G)$. If $a=1$ and $b=k$, then an $[a,b]$-factor is a $[1,k]$-factor. If $n$ is even and $k=1$, then a $[1,k]$-factor is a 1-factor (or a
perfect matching). A spanning tree of $G$ is called a spanning $k$-tree if its maximum degree is at most $k$, where $k\geq2$ is an integer.
It is obvious that a graph $G$ admits a spanning $k$-tree if and only if $G$ contains a connected $[1,k]$-factor.

Spanning $k$-trees and $[a,b]$-factors have attracted many researchers' attention. Many sufficient conditions for graphs to possess
$[1,2]$-factors were derived by Johnson, Paulusma and Wood \cite{JPW}, Kelmans \cite{Ke}, Kawarabayashi, Matsuda, Oda and Ota \cite{KMOO},
Gao and Wang \cite{GW}, Wu \cite{Wp}, Zhou \cite{Zhr,Zd}, Zhou, Zhang and Sun \cite{ZZS}, Zhou, Sun and Bian \cite{ZSB},
Zhou, Wu and Bian \cite{ZWB}, Zhou, Sun and Liu \cite{ZSL,ZSL1}, Wang and Zhang \cite{WZi}, Liu and Pan \cite{LP}. Kim, O, Park and Ree \cite{KOPR},
Wang and Zhang \cite{WZs} provided some sufficient conditions for the existence of $[1,b]$-factors in graphs. Zhou and Liu \cite{ZL} investigated
the connections between the spectral radius of a connected graph and its odd $[1,b]$-factors, and provided a lower bound on the existence of
odd $[1,b]$-factors by virtue of the spectral radius. Matsuda \cite{M}, Furuya and Yashima \cite{FY}, Wu \cite{Wa},
Zhou \cite{Za1}, Zhou, Pan and Xu \cite{ZPX2}, Zhou, Zhang and Liu \cite{ZZL}, Zhou, Pan and Xu \cite{ZPX1}, Zhou and Zhang \cite{ZZ1}, Wang and Zhang \cite{WZo,WZr} studied the properties of $[a,b]$-factors in graphs, and
showed some graphic parameter conditions for graphs to possess $[a,b]$-factors. Ding, Johnson and Seymour \cite{DJS}, Kyaw \cite{Kyaw1}, Neumann
and Eduardo \cite{NE} provided some sufficient conditions for connected graphs to possess spanning trees. Broersma and Tuinstra \cite{BT}, Tsugaki
and Yamashita \cite{TY} established the connections between degree conditions and spanning trees in connected graphs. Kyaw \cite{Kyaw2} presented
a degree and neighborhood condition for a connected graph to possess a spanning $k$-tree. Neumann-Lara and Rivera-Campo \cite{NR} gave an
independence number condition to guarantee that a connected graph contains a spanning $k$-tree. Win \cite{W} create a connection between toughness
and the existence of a spanning $k$-tree of a connected graph. Very recently, Fan, Goryainov, Huang and Lin \cite{FGHL} established a lower bound
on a spanning $k$-tree of a connected graph by virtue of its spectral radius.

Kaneko \cite{Ks} introduced the concept of leaf degree of a spanning tree. Let $T$ be a spanning tree of a connected graph $G$. The leaf degree
of a vertex $v\in V(T)$ is defined as the number of leaves adjacent to $v$ in $T$. Furthermore, the leaf degree of $T$ is the maximum leaf
degree among all the vertices of $T$. Kaneko \cite{Ks} obtained a criterion for the existence of a spanning tree with leaf degree at most $k$ in
a connected graph.

Motivated by \cite{Ks,Os} directly, it is natural and interesting to put forward a sufficient condition to guarantee that a graph contains a
spanning tree with leaf degree at most $k$. In this paper, we study the relation between the signless Laplacian spectral radius of a connected
graph and its spanning tree with leaf degree at most $k$, and obtain a signless Laplacian spectral radius condition for the existence of a spanning
tree with leaf degree at most $k$ in a connected graph, which are shown in the following.

\medskip

\noindent{\textbf{Theorem 1.1.}} Let $k$ be a positive integer, and let $G$ be a connected graph of order $n\geq k+2$.

(\romannumeral1) For $(k,n)\notin\{(1,6),(2,8),(3,10),(4,12)\}$, if $q(G)>\theta(k,n)$, then $G$ has a spanning tree with leaf degree at most
$k$, where $\theta(k,n)$ is the largest root of $x^{3}-(3n-2k-5)x^{2}+n(2n-2k-5)x-2(n-k-3)(n-k-2)=0$.

(\romannumeral2) For $(k,n)=(1,6)$, if $q(G)>4+2\sqrt{3}$, then $G$ has a spanning tree with leaf degree at most 1.

(\romannumeral3) For $(k,n)=(2,8)$, if $q(G)>5+\sqrt{21}$, then $G$ has a spanning tree with leaf degree at most 2.

(\romannumeral4) For $(k,n)=(3,10)$, if $q(G)>6+4\sqrt{2}$, then $G$ has a spanning tree with leaf degree at most 3.

(\romannumeral5) For $(k,n)=(4,12)$, if $q(G)>7+3\sqrt{5}$, then $G$ has a spanning tree with leaf degree at most 4.

\medskip

The proof of Theorem 1.1 will be given in Section 3.

\section{Some preliminaries}

In this section, we showed some necessary preliminary lemmas, which are very useful in the proof of Theorem 1.1. Kaneko \cite{Ks} provided a criterion
for the existence of a spanning tree with leaf degree at most $k$ in a connected graph.

\medskip

\noindent{\textbf{Lemma 2.1}} (\cite{Ks}). Let $k$ be a positive integer and $G$ be a connected graph. Then $G$ contains a spanning tree with leaf
degree at most $k$ if and only if
$$
i(G-S)<(k+1)|S|
$$
for every nonempty subset $S$ of $V(G)$, where $i(G-S)$ denotes the number of those isolated vertices in $G-S$.

\medskip

Cvetkovi\'c and Simi\'c \cite{CS} obtained the signless Laplacian spectral radius of a complete graph.

\medskip

\noindent{\textbf{Lemma 2.2}} (\cite{CS}). Let $K_n$ be a complete graph of order $n$, where $n\geq2$ is an integer. Then $q(K_n)=2n-2$.

\medskip

Shen, You, Zhang and Li \cite{SYZL} posed a result on the signless Laplacian spectral radius of a connected graph.

\medskip

\noindent{\textbf{Lemma 2.3}} (\cite{SYZL}). Let $G$ be a connected graph. If $H$ is a subgraph of $G$, then $q(H)\leq q(G)$. If $H$ is a
proper subgraph of $G$, then $q(H)<q(G)$.

\medskip

Let $M$ be a real symmetric matrix whose rows and columns are indexed by $V=\{1,2,\cdots,n\}$. Assume that $M$, with respect to the partition
$\pi: V=V_1\cup V_2\cup\cdots\cup V_t$, can be written as
\begin{align*}
M=\left(
  \begin{array}{ccc}
    M_{11} & \cdots & M_{1t}\\
    \vdots & \ddots & \vdots\\
    M_{t1} & \cdots & M_{tt}\\
  \end{array}
\right),
\end{align*}
where $M_{ij}$ denotes the submatrix (block) of $M$ formed by rows in $V_i$ and columns in $V_j$. Let $q_{ij}$ denote the average row sum of
$M_{ij}$. Then matrix $M_{\pi}=(q_{ij})$ is called the quotient matrix of $M$. If the row sum of every block $M_{ij}$ is a constant, then the
partition is equitable.

\medskip

\noindent{\textbf{Lemma 2.4}} (\cite{YYSX}).  Let $M$ be a real symmetric matrix with an equitable partition $\pi$, and
let $M_{\pi}$ be the corresponding quotient matrix. Then every eigenvalue of $M_{\pi}$ is an eigenvalue of $M$. Furthermore, if $M$ is
nonnegative, then the largest eigenvalues of $M$ and $M_{\pi}$ are equal.

\section{The proof of Theorem 1.1}

In this section, we prove Theorem 1.1.

\medskip

\noindent{\it Proof of Theorem 1.1.} Let $\varphi(x)=x^{3}-(3n-2k-5)x^{2}+n(2n-2k-5)x-2(n-k-3)(n-k-2)$ and let $\theta(k,n)$ is the largest root
of $\varphi(x)=0$. Suppose that $G$ has no spanning tree with leaf degree at most $k$. It follows from Lemma 2.1 that
$$
i(G-S)\geq(k+1)|S|
$$
for some nonempty subset $S$ of $V(G)$. Choose such a connected graph $G$ so that its signless Laplacian spectral radius is as large as possible.
In light of Lemma 2.3 and the choice of $G$, the induced subgraph $G[S]$ and each connected component of $G-S$ are complete graphs, respectively.
Furthermore, $G$ is just the graph $G[S]\vee(G-S)$.

We may see that $G-S$ has at most one non-trivial connected component. Otherwise, we can add edges among all non-trivial connected components to
get a non-trivial connected component of larger size, which gives a contradiction (based on Lemma 2.3). For convenience, let $i(G-S)=i$ and $|S|=s$.
Thus, we possess
\begin{align}\label{eq:3.1}
i\geq(k+1)s.
\end{align}
In what follows, we proceed by considering two possible cases.

\noindent{\bf Case 1.} $G-S$ admits just one non-trivial connected component.

In this case, it is obvious that $G=K_s\vee(K_{n_1}\cup iK_1)$, where $n_1\geq2$ is an integer. We are to show $i=(k+1)s$. Assume $i\geq(k+1)s+1$,
then we construct a new graph $H_1$ derived from $G$ by joining each vertex of $K_{n_1}$ with one vertex in $iK_1$ by an edge. Then
$i(H_1-S)\geq(k+1)s$ and $G$ is a proper spanning subgraph of $H_1$. Utilizing Lemma 2.3, we have $q(G)<q(H_1)$, which is a contradiction to
the choice of $G$. Hence, $i\leq(k+1)s$. Combining this with (\ref{eq:3.1}), we get $i=(k+1)s$, and so $G=K_s\vee(K_{n_1}\cup(k+1)sK_1)$. Obviously,
$n=(k+2)s+n_1\geq(k+2)s+2\geq k+4$. In view of the equitable partition $V(G)=V(K_s)\cup V(K_{n_1})\cup V((k+1)sK_1)$, the quotient matrix of $Q(G)$
equals
\begin{align*}
M_1=\left(
  \begin{array}{ccc}
    n+s-2 & n-(k+2)s & (k+1)s\\
    s & 2n-(2k+3)s-2 & 0\\
    s & 0 & s\\
  \end{array}
\right),
\end{align*}
whose characteristic polynomial is
\begin{align*}
f_1(x)=&x^{3}-(3n-(2k+1)s-4)x^{2}+(2n^{2}-(2k-1)sn-6n-4(k+1)s^{2}+4ks+4)x\\
&-s(2n^{2}-4(k+1)sn-6n+2(k+1)^{2}s^{2}+6(k+1)s+4).
\end{align*}
According to Lemma 2.4, the largest root, say $\theta_1$, of $f_1(x)=0$ equals $q(G)$.

Note that $K_s\vee(n-s)K_1$ is a proper subgraph of $G$. In terms of the equitable partition $V(K_s\vee(n-s)K_1)=V(K_s)\cup V((n-s)K_1)$, the
quotient matrix of $Q(K_s\vee(n-s)K_1)$ is
\begin{align*}
M_2=\left(
  \begin{array}{ccc}
    n+s-2 & n-s \\
    s & s \\
  \end{array}
\right).
\end{align*}
Then the characteristic polynomial of the matrix $M_2$ is given as
$$
f_2(x)=x^{2}-(n+2s-2)x+2s(s-1).
$$
Using Lemma 2.4, the largest root, say $\theta_2$, of $f_2(x)=0$ equals $q(K_s\vee(n-s)K_1)$. Thus, we get
\begin{align}\label{eq:3.2}
q(K_s\vee(n-s)K_1)=\theta_2=\frac{n+2s-2+\sqrt{(n+2s-2)^{2}-8s(s-1)}}{2}.
\end{align}
Together with Lemma 2.3, we infer
\begin{align}\label{eq:3.3}
\theta_1=q(G)>q(K_s\vee(n-s)K_1)=\frac{n+2s-2+\sqrt{(n+2s-2)^{2}-8s(s-1)}}{2}.
\end{align}

To verify $q(G)\leq\theta(k,n)$, it suffices to prove $\varphi(\theta_1)\leq0$. Recall that $f_1(\theta_1)=0$. Hence,
\begin{align}\label{eq:3.4}
\varphi(\theta_1)=&\varphi(\theta_1)-f_1(\theta_1)\nonumber\\
=&(s-1)(-(2k+1)\theta_1^{2}+((2k-1)n+4(k+1)s+4)\theta_1\nonumber\\
&+2(n-k-3)(n-k-2)-4(k+1)sn+2(k+1)^{2}s(s+1)+6(k+1)s)\nonumber\\
=&(s-1)g_1(\theta_1),
\end{align}
where $g_1(\theta_1)=-(2k+1)\theta_1^{2}+((2k-1)n+4(k+1)s+4)\theta_1+2(n-k-3)(n-k-2)-4(k+1)sn+2(k+1)^{2}s(s+1)+6(k+1)s$. According to (\ref{eq:3.3})
and $n\geq(k+2)s+2$, we possess
$$
\frac{(2k-1)n+4(k+1)s+4}{2(2k+1)}<n+s-2<\frac{n+2s-2+\sqrt{(n+2s-2)^{2}-8s(s-1)}}{2}<\theta_1.
$$
Thus, we derive
\begin{align}\label{eq:3.5}
g_1(\theta_1)<&g_1\left(\frac{n+2s-2+\sqrt{(n+2s-2)^{2}-8s(s-1)}}{2}\right)\nonumber\\
=&n^{2}-(4k+5)sn-(2k+5)n+2(k+1)(k+3)s^{2}+(2k^{2}+10k+10)s\nonumber\\
&+2k^{2}+6k+6-(n-s-2k-3)\sqrt{(n+2s-2)^{2}-8s(s-1)}.
\end{align}
Recall that $n\geq(k+2)s+2\geq3s+2$. For $s\geq2$, we deduce
$$
n-s-2k-3\geq(k+2)s+2-s-2k-3=(k+1)s-2k-1\geq2(k+1)-2k-1>0
$$
and
$$
(n+2s-2)^{2}-8s(s-1)-(n+s-1)^{2}=(s-1)(2n-5s-3)>0.
$$
Together with (\ref{eq:3.5}), $s\geq2$ and $n\geq(k+2)s+2$, we obtain
\begin{align}\label{eq:3.6}
g_1(\theta_1)<&n^{2}-(4k+5)sn-(2k+5)n+2(k+1)(k+3)s^{2}+(2k^{2}+10k+10)s\nonumber\\
&+2k^{2}+6k+6-(n-s-2k-3)\sqrt{(n+2s-2)^{2}-8s(s-1)}\nonumber\\
<&n^{2}-(4k+5)sn-(2k+5)n+2(k+1)(k+3)s^{2}+(2k^{2}+10k+10)s\nonumber\\
&+2k^{2}+6k+6-(n-s-2k-3)(n+s-1)\nonumber\\
=&-(4k+5)sn-n+(2k^{2}+8k+7)s^{2}+(2k^{2}+12k+12)s+2k^{2}+4k+3\nonumber\\
\leq&-(4k+5)s((k+2)s+2)-((k+2)s+2)+(2k^{2}+8k+7)s^{2}\nonumber\\
&+(2k^{2}+12k+12)s+2k^{2}+4k+3\nonumber\\
=&-(2k^{2}+5k+3)s^{2}+(2k^{2}+3k)s+2k^{2}+4k+1\nonumber\\
\leq&-2(2k^{2}+5k+3)s+(2k^{2}+3k)s+2k^{2}+4k+1\nonumber\\
=&-(2k^{2}+7k+6)s+2k^{2}+4k+1\nonumber\\
\leq&-2(2k^{2}+7k+6)+2k^{2}+4k+1\nonumber\\
=&-2k^{2}-10k-11\nonumber\\
<&0.
\end{align}
It follows from (\ref{eq:3.4}), (\ref{eq:3.6}) and $s\geq1$ that
$$
\varphi(\theta_1)=(s-1)g_1(\theta_1)\leq0,
$$
which yields $q(G)=\theta_1\leq\theta(k,n)$, which contradicts $q(G)>\theta(k,n)$ for $n\geq k+2$ and $(k,n)\notin\{(1,6),(2,8),(3,10),(4,12)\}$.
As for $(k,n)=(1,6)$, we possess $q(G)=\theta_1\leq\theta(1,6)<4+2\sqrt{3}$, which is a contradiction to $q(G)>4+2\sqrt{3}$ for $(k,n)=(1,6)$.
As for $(k,n)=(2,8)$, we have $q(G)=\theta_1\leq\theta(2,8)<5+\sqrt{21}$, which contradicts $q(G)>5+\sqrt{21}$ for $(k,n)=(2,8)$. As for
$(k,n)=(3,10)$, we obtain $q(G)=\theta_1\leq\theta(3,10)<6+4\sqrt{2}$, which contradicts $q(G)>6+4\sqrt{2}$ for $(k,n)=(3,10)$. As for $(k,n)=(4,12)$,
we derive $q(G)=\theta_1\leq\theta(4,12)<7+3\sqrt{5}$, which is a contradiction to $q(G)>7+3\sqrt{5}$ for $(k,n)=(4,12)$.

\noindent{\bf Case 2.} $G-S$ has no non-trivial connected component.

In this case, it is clear that $G=K_s\vee iK_1$. We are to show $i\leq(k+1)s+1$. Assume that $i\geq(k+1)s+2$. Then we can construct a new graph
$H_2$ obtained from $G$ by adding an edge in $iK_1$. Thus, we see $i(H_2-S)\geq(k+1)s$ and $G$ is a proper subgraph of $H_2$. Using Lemma 2.3,
we have $q(G)<q(H_2)$, which contradicts the choice of $G$. Therefore, $i\leq(k+1)s+1$. Together with (\ref{eq:3.1}), it suffices to consider
$i=(k+1)s$ and $i=(k+1)s+1$.

\noindent{\bf Subcase 2.1.} $i=(k+1)s$.

In this subcase, we have $G=K_s\vee(k+1)sK_1=K_s\vee(n-s)K_1$, where $n=(k+2)s\geq k+2$. In terms of (\ref{eq:3.2}), we get
$$
q(G)=q(K_s\vee(n-s)K_1)=\theta_2=\frac{n+2s-2+\sqrt{(n+2s-2)^{2}-8s(s-1)}}{2}.
$$
Note that $f_2(\theta_2)=0$. Then
\begin{align}\label{eq:3.7}
\varphi(\theta_2)=&\varphi(\theta_2)-\theta_2f_2(\theta_2)\nonumber\\
=&-(2n-2s-2k-3)\theta_2^{2}+(2n^{2}-2kn-5n-2s^{2}+2s)\theta_2-2(n-k-3)(n-k-2)\nonumber\\
=&-(s+1)n^{2}+(3s+2k+2)sn+(2k+5)n-(2s^{2}+4k+4)s-2k^{2}-6k-6\nonumber\\
&+(-sn+n+s^{2}+2ks+2s-2k-3)\sqrt{(n+2s-2)^{2}-8s(s-1)}.
\end{align}

\noindent{\bf Subcase 2.1.1.} $s=1$.

Obviously, $n=k+2$ and $q(G)=n=\theta(k,k+2)$, which contradicts $q(G)>\theta(k,n)$ for $n=k+2$.

\noindent{\bf Subcase 2.1.2.} $s=2$.

Note that $n=2k+4$. If $k=1$, then $n=6$ and $q(G)=4+2\sqrt{3}$, a contradiction. If $k=2$, then $n=8$ and $q(G)=5+\sqrt{21}$, a contradiction.
If $k=3$, then $n=10$ and $q(G)=6+4\sqrt{2}$, a contradiction. If $k=4$, then $n=12$ and $q(G)=7+3\sqrt{5}$, a contradiction. If $k\geq5$, then
it follows from $n=2k+4$ and (\ref{eq:3.7}) that $\varphi(\theta_2)=-2k^{2}+4k+6+2\sqrt{(k+1)(k+5)}<0$, which implies that $q(G)=\theta_2<\theta(k,n)$
for $s=2$, a contradiction.

\noindent{\bf Subcase 2.1.3.} $s\geq3$.

Note that $n=(k+2)s$. We easily deduce
\begin{align*}
&-sn+n+s^{2}+2ks+2s-2k-3\\
=&-(k+2)s^{2}+(k+2)s+s^{2}+2ks+2s-2k-3\\
=&-(k+1)s^{2}+(3k+4)s-2k-3\\
<&0
\end{align*}
and
$$
(n+2s-2)^{2}-8s(s-1)-(n+s-1)^{2}=(s-1)(2n-5s-3)\geq0
$$
for $s\geq3$. Together with (\ref{eq:3.7}) and $n=(k+2)s$, we get
\begin{align}\label{eq:3.8}
\varphi(\theta_2)\leq&-(s+1)n^{2}+(3s+2k+2)sn+(2k+5)n-(2s^{2}+4k+4)s-2k^{2}-6k-6\nonumber\\
&+(-sn+n+s^{2}+2ks+2s-2k-3)(n+s-1)\nonumber\\
=&-(2k^{2}+5k+3)s^{3}+(4k^{2}+16k+13)s^{2}-7(k+1)s-2k^{2}-4k-3\nonumber\\
:=&p(s).
\end{align}
Then we see that $p'(s)=-3(2k^{2}+5k+3)s^{2}+2(4k^{2}+16k+13)s-7(k+1)$ and $p''(s)=-6(2k^{2}+5k+3)s+2(4k^{2}+16k+13)$. For $s\geq3$, we infer
\begin{align*}
p''(s)=&-6(2k^{2}+5k+3)s+2(4k^{2}+16k+13)\\
\leq&-18(2k^{2}+5k+3)+2(4k^{2}+16k+13)\\
=&-28k^{2}-58k-28\\
<&0.
\end{align*}
Hence, $p'(s)$ is decreasing in the interval $[3,+\infty)$. Then we have
\begin{align*}
p'(s)\leq&p'(3)\\
=&-27(2k^{2}+5k+3)+6(4k^{2}+16k+13)-7(k+1)\\
=&-30k^{2}-46k-10\\
<&0,
\end{align*}
which implies that $p(s)$ is decreasing in the interval $[3,+\infty)$. Thus, we get
\begin{align*}
p(s)\leq&p(3)\\
=&-27(2k^{2}+5k+3)+9(4k^{2}+16k+13)-21(k+1)-2k^{2}-4k-3\\
=&-20k^{2}-16k+12\\
<&0.
\end{align*}
Together with (\ref{eq:3.8}), we obtain $\varphi(\theta_2)\leq p(s)<0$ for $s\geq3$, which implies $q(G)=\theta_2<\theta(k,n)$, a contradiction.

\noindent{\bf Subcase 2.2.} $i=(k+1)s+1$.

In this subcase, it is obvious that $G=K_s\vee((k+1)s+1)K_1$ and $n=(k+2)s+1\geq k+3$. According to the equitable partition $V(G)=V(K_s)\cup V(((k+1)s+1)K_1)$,
the quotient matrix of $Q(G)$ equals
\begin{align*}
M_3=\left(
  \begin{array}{ccc}
    (k+3)s-1 & (k+1)s+1\\
    s & s \\
  \end{array}
\right).
\end{align*}
Thus we derive its characteristic polynomial as
$$
f_3(x)=x^{2}-((k+4)s-1)x+2s(s-1).
$$
By virtue of Lemma 2.4, the largest root, say $\theta_3$, of $f_3(x)=0$ equals $q(G)$. Thus, we possess
$$
q(G)=\theta_3=\frac{(k+4)s-1+\sqrt{(k^{2}+8k+8)s^{2}-2ks+1}}{2}.
$$
Note that $n=(k+2)s+1$ and $f_3(\theta_3)=0$. Then
\begin{align}\label{eq:3.9}
\varphi(\theta_3)=&\varphi(\theta_3)-\theta_3f_3(\theta_3)\nonumber\\
=&-(2(k+1)s-2k-1)\theta_3^{2}+(2(k+1)(k+3)s^{2}-k(2k+5)s-2k-3)\theta_3\nonumber\\
&-2(k+2)^{2}s^{2}+2(k+2)(2k+3)s-2(k+1)(k+2)\nonumber\\
=&-k(k+1)s^{3}+(k^{2}-1)s^{2}+(2k^{2}+5k+3)s-2(k+1)^{2}\nonumber\\
&-(k+1)(s-1)(s-2)\sqrt{(k^{2}+8k+8)s^{2}-2ks+1}\nonumber\\
\leq&-k(k+1)s^{3}+(k^{2}-1)s^{2}+(2k^{2}+5k+3)s-2(k+1)^{2}\nonumber\\
&-(k+1)(s-1)(s-2)(k+3)s.
\end{align}

If $s=1$, then $n=k+3$. Together with (\ref{eq:3.9}), we get $q(G)=\theta_3=\theta(k,k+3)$, which contradicts $q(G)>\theta(k,n)$ for $n=k+3$.
If $s=2$, then $n=2k+5$. Combining this with (\ref{eq:3.9}), we infer
\begin{align*}
\varphi(\theta_3)=&-8k(k+1)+4(k^{2}-1)+2(2k^{2}+5k+3)-2(k+1)^{2}\\
=&-2k^{2}-2k\\
<&0,
\end{align*}
which implies $q(G)=\theta_3<\theta(k,2k+5)$, which contradicts $q(G)>\theta(k,n)$ for $n=2k+5$. If $s\geq3$, then it follows from (\ref{eq:3.9})
that
\begin{align*}
\varphi(\theta_3)\leq&-k(k+1)s^{3}+(k^{2}-1)s^{2}+(2k^{2}+5k+3)s-2(k+1)^{2}\\
&-(k+1)(s-1)(s-2)(k+3)s\\
<&-3k(k+1)s^{2}+(k^{2}-1)s^{2}+(2k^{2}+5k+3)s-2(k+1)^{2}\\
=&-(2k^{2}+3k+1)s^{2}+(2k^{2}+5k+3)s-2(k+1)^{2}\\
\leq&-3(2k^{2}+3k+1)s+(2k^{2}+5k+3)s-2(k+1)^{2}\\
=&-4k(k+1)s-2(k+1)^{2}\\
<&0,
\end{align*}
which yields $q(G)=\theta_3<\theta(k,n)$, a contradiction. This completes the proof of Theorem 1.1. \hfill $\Box$

\medskip

\section{Extremal graphs}

Finally, we establish several graphs to claim that the bounds obtained in Theorem 1.1 are sharp.

\medskip

\noindent{\textbf{Theorem 4.1.}} Let $k$ and $n$ be two positive integers with $n\geq k+2$, and let $\theta(k,n)$ be the largest root of
$x^{3}-(3n-2k-5)x^{2}+n(2n-2k-5)x-2(n-k-3)(n-k-2)=0$. If $(k,n)\notin\{(1,6),(2,8),(3,10),(4,12)\}$, then
$q(K_1\vee(K_{n-k-2}\cup(k+1)K_1))=\theta(k,n)$ and $K_1\vee(K_{n-k-2}\cup(k+1)K_1)$ contains no spanning tree with leaf degree at most $k$. If
$(k,n)=(1,6)$, then $q(K_2\vee4K_1)=4+2\sqrt{3}$ and $K_2\vee4K_1$ contains no spanning tree with leaf degree at most 1. If $(k,n)=(2,8)$, then
$q(K_2\vee6K_1)=5+\sqrt{21}$ and $K_2\vee6K_1$ contains no spanning tree with leaf degree at most 2. If $(k,n)=(3,10)$, then $q(K_2\vee8K_1)=6+4\sqrt{2}$
and $K_2\vee8K_1$ contains no spanning tree with leaf degree at most 3. If $(k,n)=(4,12)$, then $q(K_2\vee10K_1)=7+3\sqrt{5}$ and $K_2\vee10K_1$
contains no spanning tree with leaf degree at most 4.

\medskip

\noindent{\it Proof.} In terms of the equitable partition $V(K_1\vee(K_{n-k-2}\cup(k+1)K_1))=V(K_1)\cup V(K_{n-k-2})\cup V((k+1)K_1))$, the quotient
matrix of $Q(K_1\vee(K_{n-k-2}\cup(k+1)K_1))$ equals
\begin{align*}
M(K_1\vee(K_{n-k-2}\cup(k+1)K_1))=\left(
  \begin{array}{ccc}
    n-1 & n-k-2 & k+1\\
    1 & 2n-2k-5 & 0\\
    1 & 0 & 1\\
  \end{array}
\right),
\end{align*}
whose characteristic polynomial is given as $x^{3}-(3n-2k-5)x^{2}+n(2n-2k-5)x-2(n-k-3)(n-k-2)$. In view of Lemma 2.4, the largest root of $x^{3}-(3n-2k-5)x^{2}+n(2n-2k-5)x-2(n-k-3)(n-k-2)=0$, $\theta(k,n)$ equals $q(K_1\vee(K_{n-k-2}\cup(k+1)K_1))$. Namely,
$q(K_1\vee(K_{n-k-2}\cup(k+1)K_1))=\theta(k,n)$. Set $S=V(K_1)$. Then $|S|=1$ and
$i(K_1\vee(K_{n-k-2}\cup(k+1)K_1)-S)=i((K_{n-k-2}\cup(k+1)K_1))=k+1=(k+1)|S|$. According to Lemma 2.1, $K_1\vee(K_{n-k-2}\cup(k+1)K_1)$ has no
spanning tree with leaf degree at most $k$.

\medskip

By virtue of the equitable partition $V(K_2\vee4K_1)=V(K_2)\cup V(4K_1)$, the quotient matrix of $Q(K_2\vee4K_1)$ equals
\begin{align*}
M(K_2\vee4K_1)=\left(
  \begin{array}{ccc}
    6 & 4\\
    2 & 2 \\
  \end{array}
\right).
\end{align*}
Then its characteristic polynomial is $x^{2}-8x+4$. By means of Lemma 2.4, the largest root of $x^{2}-8x+4=0$ equals $q(K_2\vee4K_1)$. Namely,
$q(K_2\vee4K_1)=4+2\sqrt{3}$. Write $S=V(K_2)$. Then $|S|=2$ and $i(K_2\vee4K_1-S)=4=2|S|$. Applying Lemma 2.1, $K_2\vee4K_1$ has no spanning
tree with leaf degree at most 1.

\medskip

According to the equitable partition $V(K_2\vee6K_1)=V(K_2)\cup V(6K_1)$, the quotient matrix of $Q(K_2\vee6K_1)$ equals
\begin{align*}
M(K_2\vee6K_1)=\left(
  \begin{array}{ccc}
    8 & 6\\
    2 & 2 \\
  \end{array}
\right).
\end{align*}
Then its characteristic polynomial is $x^{2}-10x+4$. Using Lemma 2.4, the largest root of $x^{2}-10x+4=0$ equals $q(K_2\vee6K_1)$. That is,
$q(K_2\vee6K_1)=5+\sqrt{21}$. Let $S=V(K_2)$. Then $|S|=2$ and $i(K_2\vee6K_1-S)=6>4=2|S|$. In light of Lemma 2.1, $K_2\vee6K_1$ has no spanning
tree with leaf degree at most 2.

\medskip

By the equitable partition $V(K_2\vee8K_1)=V(K_2)\cup V(8K_1)$, the quotient matrix of $Q(K_2\vee8K_1)$ is
\begin{align*}
M(K_2\vee8K_1)=\left(
  \begin{array}{ccc}
    10 & 8\\
    2 & 2 \\
  \end{array}
\right),
\end{align*}
whose characteristic polynomial is given as $x^{2}-12x+4$. According to Lemma 2.4, the largest root of $x^{2}-12x+4=0$ equals $q(K_2\vee8K_1)$.
That is, $q(K_2\vee8K_1)=6+4\sqrt{2}$. Let $S=V(K_2)$. Then $|S|=2$ and $i(K_2\vee8K_1-S)=8>4=2|S|$. Utilizing Lemma 2.1, $K_2\vee8K_1$ has no
spanning tree with leaf degree at most 3.

\medskip

In view of the equitable partition $V(K_2\vee10K_1)=V(K_2)\cup V(10K_1)$, the quotient matrix of $Q(K_2\vee10K_1)$ is
\begin{align*}
M(K_2\vee10K_1)=\left(
  \begin{array}{ccc}
    12 & 10\\
    2 & 2 \\
  \end{array}
\right),
\end{align*}
whose characteristic polynomial is given as $x^{2}-14x+4$. By Lemma 2.4, the largest root of $x^{2}-14x+4=0$ equals $q(K_2\vee10K_1)$.
That is, $q(K_2\vee10K_1)=7+3\sqrt{5}$. Put $S=V(K_2)$. Then $|S|=2$ and $i(K_2\vee10K_1-S)=10>4=2|S|$. By means of Lemma 2.1, $K_2\vee10K_1$ has
no spanning tree with leaf degree at most 4. \hfill $\Box$

\section*{Declaration of competing interest}

\medskip

The authors declares that they have no known competing financial interests or personal relationships that could have appeared to influence the work
reported in this paper.

\section*{Data availability}

\medskip

No data was used for the research described in the article.

\medskip

%\section*{Acknowledgments}

\end{document}